\documentclass[10pt,english]{article} 

\usepackage{amsmath, amsthm, amssymb}
\usepackage{a4, babel, graphicx, listings, fancyhdr, verbatim}
\usepackage{mathrsfs} 
\usepackage{anysize} 
\usepackage{natbib} 
\usepackage[format = hang]{caption}                              

\usepackage{a4}
\usepackage[latin1]{inputenc}
\usepackage{amsmath, amsthm}
\usepackage{graphicx}
\usepackage{amssymb}
\usepackage{verbatim}
\usepackage{listings}
\usepackage{alltt}
\usepackage{babel}
\usepackage{mathrsfs}
\usepackage{ae}
\usepackage{natbib}
\usepackage{booktabs} 

\newtheoremstyle{theorem}
{10pt} 
{10pt} 
{\sl} 
{\parindent} 
{\bf} 
{. } 
{ } 
{} 

\theoremstyle{theorem}

\def\beq{\begin{eqnarray}}
\def\eeq{\end{eqnarray}}

\def\beqn{\begin{eqnarray*}}  
\def\eeqn{\end{eqnarray*}}

\def\N{{\rm N}}

\def\Pr{P}

\def\hatt{\widehat}

\def\midd{\,|\,}

\def\cc{{\rm cc}}

\numberwithin{equation}{section} 
\numberwithin{figure}{section}
\numberwithin{table}{section}
\usepackage{hyperref}

\title{Confidence in confidence distributions!}

\date{February 2020}

\begin{document}


\maketitle

\centerline{\large\bf C\'eline Cunen$^1$, Nils Lid Hjort$^1$, 
   Tore Schweder$^2$}

\medskip 
\centerline{\bf $^1$Department of Mathematics, University of Oslo}

\medskip 
\centerline{\bf $^2$Department of Economics, University of Oslo}

\begin{abstract}
\noindent
The recent article `Satellite conjunction analysis and the false 
confidence theorem' (Balch, Martin, and Ferson, 2019, this journal) 
points to certain difficulties with Bayesian analysis when
used for models for satellite conjuntion and ensuing operative
decisions. Here we supplement these previous analyses and
findings with further insights, uncovering what we perceive
of as being the crucial points, explained in a prototype
setup where exact analysis is attainable. We also show
that a different and frequentist method, involving confidence
distributions, is free of the false confidence syndrome. 
\end{abstract}




\bigskip\noindent
Bayesian analysis is often employed to calculate the collision 
probability of satellites based on noisy measurement of their 
position and speed. 
Balch, Martin, Ferson (2019), below referred to as BMF, 
have demonstrated that these analyses can be highly misleading, 
and they argue that this is 
an instance of a more general phenomenon, 
which they call {\it false confidence}.
BMF find that the confidence in the event 
of no collision tends to get larger with noisier data 
even when the satellites actually are on a collision course.
Such unfortunate properties are known to be present 
in some Bayesian analyses, 
cf.~Schweder and Hjort (2016, e.g.~Chapters 14, 15). 
We argue that such probability dilution is avoided 
when frequentist analysis by way of confidence 
distributions is employed. We demonstrate this below 
by revisiting the satellite collision problem.
The purpose of this note is also to make clearer 
what `the heart of the matter' is: in certain setups,
even a sound-looking Bayesian analysis might 
have seriously unfortunate frequentist properties. 
Such consequences are particularly 
drastic when the methods are put to repeated use,
as for monitoring of satellite collision probabilities. 



Confidence distributions (CDs) are a type of inferential 
summary whose outward appearance is similar to 
a Bayesian posterior distribution. Just like a 
posterior for some parameter of interest $\delta$,
after having observed data $y$,  
a CD $C(\delta\midd y)$ defines 
a probability distribution on the space of possible 
$\delta$ values. The function $C(\delta\midd y)$ 
is a data-dependent cumulative distribution function. 
Unlike a posterior, however, a CD is required 
to have correct frequentist coverage, 
i.e.~$C(\delta_0\midd y)$ is uniformly distributed when 
data are regarded as stochastic and $\delta_0$ is 
the true value of the parameter. 
The full confidence curve $C(\delta\midd y)$ 
is a practical summary of the uncertainty associated 
with the estimated parameter; for more on CDs see 
Schweder and Hjort (2016), Hjort and Schweder (2018). 
CDs and the related concept confidence curves 
are strongly related to confidence intervals. 
Since properly calibrated confidence intervals 
have guaranteed coverage, they are free of false confidence,  
and so are CDs, 
at least with respect to all interval statements, 
including confidence intervals and p-values. 


After these brief general points we proceed to 
present an analysis of a simplified version 
of BMF's satellite example. 
This simplification allows explicit expressions 
of the functions involved, but retains what we 
see as the heart of the matter, a crucial 
and problematic aspect of the Bayesian solution. 
Note that sophisticated satellite conjunction risk assessors 
are aware of the false confidence phenomenon, 
and the related issue of probability dilution, 
which affects the Bayesian calculation of collision 
probability. In practice, they therefore use 
this calculation not as a guarantee of safety 
but rather as a way to identify conjunctions that 
are likely to be dangerous. Their null 
hypothesis is that a collision will not occur 
and when the estimated collision probability 
is high the conjunction is identified as `dangerous'. 
See Hejduk, Snow, Newman (2019) 
for an extended treatment of these issues. 

Following Martin (2019), 
assuming that the two satellites have constant 
and identical speed 
we can reduce the problem to two dimensions, 
with $\theta=(\theta_1,\theta_2)$ being the true differences 
between the unknown positions of the satellites 
along each axis on a plane. The interest lies 
in inference on the distance between the satellites, 
i.e.~$\delta=||\theta|| = (\theta_1^2 + \theta_2^2)^{1/2}$. 
We observe a single pair $y=(y_1,y_2)$ and 
assume that these two are independent and normal 
with means $(\theta_1,\theta_2)$ and equal 
known variances $\sigma^2$.  

BFM present a Bayesian framework with flat priors,  
and in this simplified setup we get the simple
posterior distribution where $\theta_i\midd y_i$
is a normal $\N(y_i,\sigma^2)$ for the two components.
From this posterior one can obtain a posterior for 
$\delta$ by integration, which is most conveniently 
presented in the form of a cumulative distribution 
function on the space of $\delta$,
\beq
\label{eq:BB}
B(\delta\midd y_1,y_2) 
   = \Gamma_2\Bigl({\delta^2\over \sigma^2},
   {||y||^2\over \sigma^2}\Bigr), 
\eeq
writing $\Gamma_2(\cdot,\nu)$ for the cumulative distribution
function of a non-central $\chi^2$ with two degrees of freedom 
and noncentrality parameter $\nu$. 
The component posterior noted above is also the natural 
CD for $\theta_1,\theta_2$ separately, but one cannot 
derive the confidence distribution of $\delta$ by 
manipulating the joint CD based on these two independent CDs.
Instead, the CD of $\delta$ is found by considering 
the sampling distribution of 
$\hatt\delta = (y_1^2 + y_2^2)^{1/2}$, 
which is a sufficient statistic for $\delta$. One finds 
\beq
\label{eq:CC}
C(\delta\midd y_1,y_2) = 1-\Gamma_2\Bigl({||y||^2\over \sigma^2},
   {\delta^2\over \sigma^2}\Bigr).
\eeq 
Incidentally, this agrees with the belief 
in $[0,\delta]$, found by the Inferential Model
methods of Martin and Liu (2015), for this situation.
Note also that the marginalisation to $\delta$ is carried out
in the the observation space and not in the epistemic
probability space 
(which is what is done in the Bayesian analysis). 

\begin{figure}[h]
\centering
\includegraphics[scale=0.66]{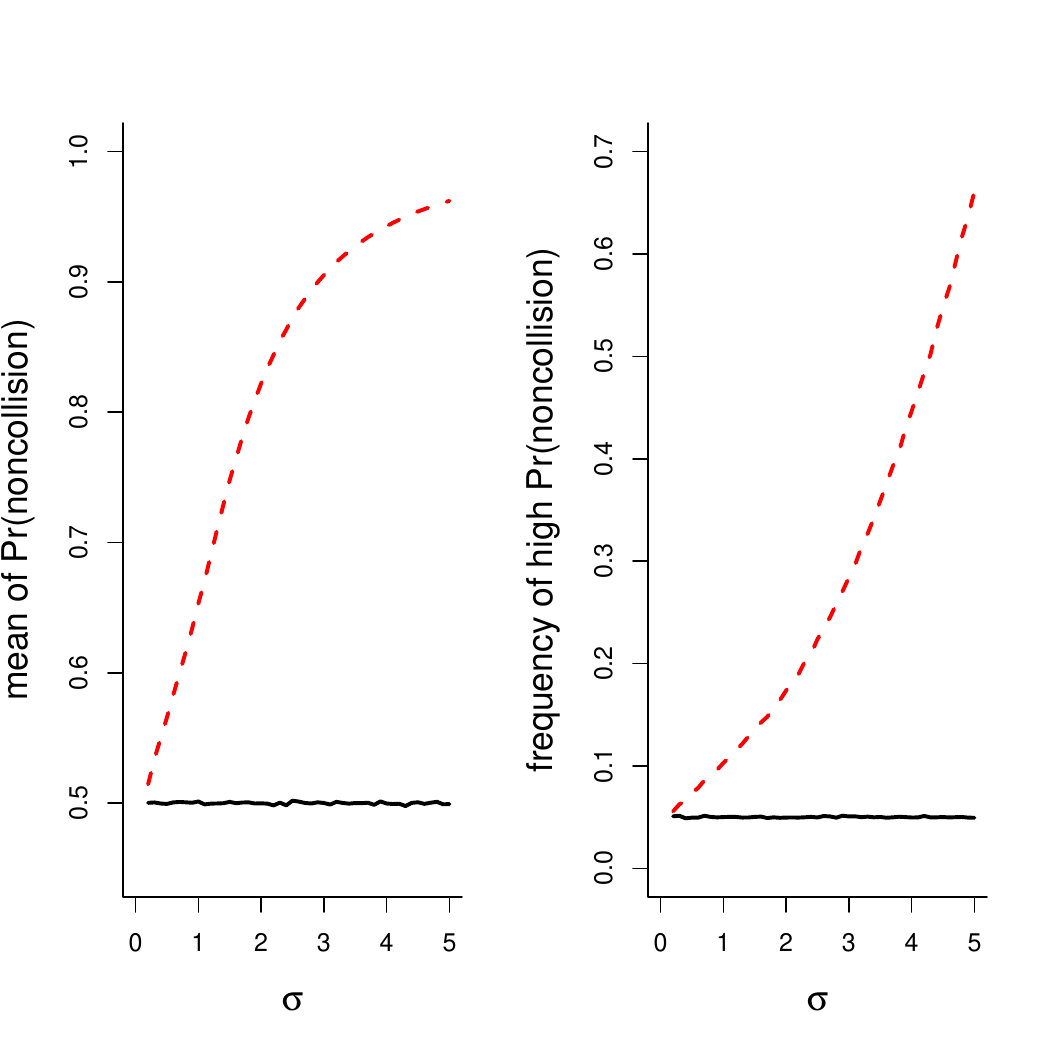}
\caption{
For each value of $\sigma$, we have computed 
the distribution of non-collision probabilities 
$1-B(2.00\midd y_1,y_2)$ and $1-C(2.00\midd y_1,y_2)$,
with the Bayesian (slanted, red curves) 
and frequentist CD (full, black curves) methods. 
This is for a setup with $R=2.00$ the threshold 
for collision and true value $\delta=1.99$.  
Left panel: the means of these probabilities; 
right panel: the frequency of high probabilities,
those above 0.95.} 
\label{figure:fig12}
\end{figure}

Let us briefly investigate the non-collision probabilities 
that these two methods produce.
Collision is defined as the event that $\delta$ is smaller 
than $R$, the combined radius of the two satellites.  
We set $R=2.00$ and $\delta=1.99$ -- a setup where 
the satellites really are on a collision 
course, but just barely so. For each of a set of 
values of $\sigma$ we simulate $10^5$ 
realisations of normal pairs $(y_1,y_2)$.
For each such dataset we then calculate the probability 
of non-collision, according to each method, 
i.e.~$\Pr(\delta > 2.00\midd y_1,y_2)$, 
which is equal to $1-B(2.00\midd y_1,y_2)$ for the Bayesian method, 
and $1-C(2.00\midd y_1,y_2)$ for the CD. 
Figure \ref{figure:fig12} displays 
frequentist properties of the Bayesian 
and the frequentist CD curves for two aspects 
of these non-collision probabilities; 
the left panel shows their mean non-collision
probabilities, the right panel the frequency of high values,
above 0.95. Exact calculations are incidentally 
also possible here. 


The misleading behaviour of the Bayesian solution 
is apparent: as the observation noise increases, 
the Bayesian method reports higher and higher 
confidence in non-collision. In particular, the Bayesian
method is {\it biased} in repeated use. 
The CD, on the other hand, reports probabilities that are 
correctly calibrated, in the sense that it 
wrongly indicates a high confidence in non-collision 
for only 5\% of the datasets. A mean probability 
of non-collision of around 50\% might seem high, 
but keep in mind that the data are drawn from 
a setup where the true distance is such that 
the satellites barely collide. When the observation 
error is high there is a high chance of observing 
data that (wrongly) indicate a large estimated distance.
This misleading property of the Bayesian probability 
of the satellites to not collide was identified by BMF 
and called false confidence.

The frequentist confidence $1-C(R\midd y_1,y_2)$ 
is free of such false confidence. More formally we actually 
have that $1-C(R\midd y_1,y_2)$ 
is uniformly distributed if the true minimal distance is $R$, 
i.e.~if $\delta=R$. 
If however $\delta<R$ the distribution of $1-C(R\midd y_1,y_2)$ 
will be shifted to the left of the uniform, 
and towards larger values if $\delta>R$. 
Note furthermore that $C(R\midd y_1,y_2)$ is the confidence 
we have in $\delta\in[0,R]$, i.e.~in collision, 
after having observed our data; also, $1-C(R\midd y_1,y_2)$
is the p-value for testing the null hypothesis that 
there will be no collision (see Schweder, 2018). 
Again the frequentist CD is fully calibrated, whereas the Bayesian machine 
provides misleading answers, the more so for increasing 
noise level. 

\begin{figure}[h]
\centering
\includegraphics[scale=0.66]{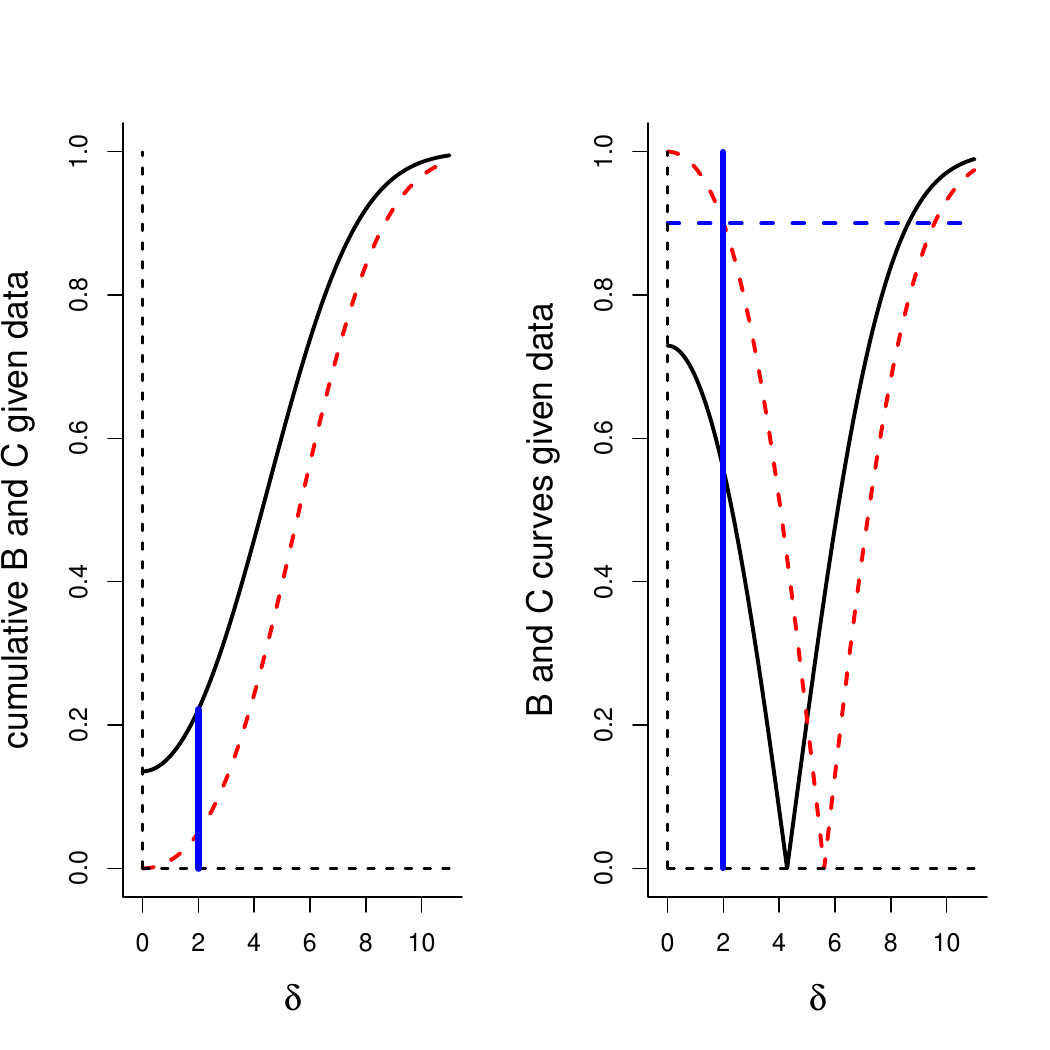}
\caption{
Left panel: the CD (black curve) and the Bayesian
posterior cumulative for $\delta$, after having 
observed $\|y\|=5.00$, with assumed $\sigma=2.50$;
the critical value is $R=2.00$ (marked blue),
where we read off the confidence 0.222 in $[0,R]$.  
Right panel: the corresponding confidence curve 
$\cc(\delta\midd y_1,y_2)=|1-2\,C(\delta\midd y_1,y_2)|$
(black) and the Bayesian credibility curve 
$|1-2\,B(\delta\midd y_1,y_2)|$ (red). 
Reading off 90\% confidence and credibility intervals 
yields $[0.00,8.63]$ and $[2.01,9.57]$, respectively.
The true $\delta=1.99$ behind the generation of 
$\|y\|$ here is indicated by the vertical blue line.}
\label{figure:fig13}
\end{figure}

It is also illuminating to see how the Bayesian cumulative  
$B(\delta\midd y_1,y_2)$ and frequentist $C(\delta\midd y_1,y_2)$ 
of (\ref{eq:BB}) and (\ref{eq:CC}) pan out in practice.
In Figure \ref{figure:fig13} we display these cumulatives
(left panel), in a situation with observed length $\|y\|=5.00$,
with assumed $\sigma=2.50$. The right panel then 
shows the useful Bayesian and confidence curves, 
$|1-2\,B(\delta\midd y_1,y_2)|$ and $|1-2\,C(\delta\midd y_1,y_2)|$.
These are convenient data summaries, for the most 
pertinent parameter, the $\delta$. They `point' to 
the median confidence estimates, 4.29 for the CD and 
5.61 for Bayes, and also make it easy to read off 
confidence intervals. Here, with coverage level 0.90, 
the Bayesian interval $[2.01,9.57]$ is unfortunate, 
in that it misses the true value $\delta=1.99$; 
the $\cc(\delta\midd y_1,y_2)$ fares rather better 
with its $[0.00,8.63]$ interval. 

The frequentist properties, also for a Bayesian statistical 
method, are of interest in contexts like the satellite 
collision problem, where potential collision events occur frequently.
When the model parameter is given a vague prior, as there, 
but the basic parameter of interest is a nonlinear 
function of this underlying model parameter, 
problems like false confidence and bias might be present.
When  confidence can be measured in the frequentist 
spirit of Fisher and Neyman, as laid out in 
Schweder and Hjort (2016), there will never 
be any false confidence, and we can trust the obtained confidence!

\medskip 
{\bf Acknowledgements.}
We are grateful for comments from both Ryan Martin 
and from anonymous referees, which have contributed 
to a clearer presentation. 

\bigskip
\centerline{\bf References} 

\def\ref#1{{\noindent\hangafter=1\hangindent=20pt
  #1\smallskip}}          
\parindent0pt
\baselineskip11pt
\parskip3pt 

\medskip
\ref{%
Balch, M.S., Martin, R., and Ferson, S. (2019). 
Satellite conjunction analysis and the false confidence theorem.
{\sl Proceedings of the Royal Society, A}, 
475, issue 2227.}

\ref{%
Hejduk, M.D., Snow, D.E., and Newman, L.K. (2019). 
Satellite conjunction assessment risk analysis for 
`dilution region' events: issues and operational approaches. 
In {\sl Space Traffic Management Conference} {\bf 28}, 
Austin, TX, February. 
{\tt commons.erau.edu/cgi/viewcontent.cgi?article=1294\allowbreak\&context=stm} }

\ref{%
Hjort, N.L.~and Schweder, T. (2018). 
Confidence distributions and related themes. 
[General introduction article to a Special Issue, 
dedicated to this topic.] 
{\sl Journal of Statistical Planning and Inference}
{\bf 195}, 1--13.}

\ref{%
Martin, R. (2019). 
False confidence, non-additive beliefs, and valid statistical inference. 
{\sl International Journal of Approximate Reasoning}, 113, 39--73.}

\ref{%
Martin, R.~and Liu, C. (2015).
{\sl Inferential Models: Reasoning with Uncertainty.}
CRS Press, New York.} 

\ref{%
Schweder, T. (2018).
Confidence is epistemic probability for empirical science.
{\sl Journal of Statistical Planning and Inference} {\bf 195},
116--125.} 

\ref{%
Schweder, T.~and Hjort, N.L. (2016).
{\sl Confidence, Likelihood, Probability: 
Statistical Inference with Confidence Distributions.} 
Cambridge University Press.}

\end{document}